\documentclass[12pt]{article}
\usepackage{amsmath}
\usepackage{latexsym}
\usepackage{amssymb}
\usepackage{amsfonts}
\usepackage{amsthm}
\title{PURE DIFFERENTIAL MODULES \\ AND A RESULT OF MACAULAY \\  ON UNMIXED POLYNOMIAL IDEALS}
\author{J.-F. Pommaret \\ CERMICS, Ecole des Ponts ParisTech,\\ 6/8 Av. Blaise Pascal, 77455 Marne-la-Vall\'ee Cedex 02, France \\
E-mail: jean-francois.pommaret@wanadoo.fr, pommaret@cermics.enpc.fr \\
URL: http://cermics.enpc.fr/$\sim$pommaret/home.html }
\date{  }
\begin{document}
\maketitle

\noindent
{\bf ABSTRACT} \\ 
 
The first purpose of this paper is to point out a curious result announced by Macaulay on the Hilbert function of a differential module in his famous book {\it  The Algebraic Theory of Modular Systems } published in 1916. Indeed, on page $78/79$ of this book, Macaulay is saying the following:\\

" A polynomial ideal $\mathfrak{a} \subset k[{\chi}_1$,..., ${\chi}_n]=k[\chi]$ is of the {\it principal class} and thus {\it unmixed} if it has rank $r$ and is generated by $r$ polynomials. Having in mind this definition, a primary ideal $\mathfrak{q}$ with associated prime ideal $\mathfrak{p} = rad(\mathfrak{q})$ is such that any ideal $\mathfrak{a}$ of the principal class with $\mathfrak{a} \subset \mathfrak{q}$ determines a primary ideal of greater {\it multiplicity} over $k$. In particular, we have $dim_k(k[\chi]/({\chi}_1$,...,${\chi}_n)^2)=n+1$ because, passing to a system of PD equations for one unknown $y$, the parametric jets are \{${y,y_1, ... ,y_n}$\} but any ideal $\mathfrak{a}$ of the principal class with $\mathfrak{a}\subset ({\chi}_1,{É},{\chi}_n)^2$ is contained into a {\it simple} ideal, that is a primary ideal 
$\mathfrak{q}$ such that $rad(\mathfrak{q})=\mathfrak{m}\in max(k[\chi])$ is a maximal and thus prime ideal with $dim_k(M)=dim_k(k[\chi]/\mathfrak{q})=2^n$ at least. Accordingly, any primary ideal $\mathfrak{q}$ may not be a member of the primary decomposition of an unmixed ideal $\mathfrak{a} \subseteq \mathfrak{q}$ of the principal class. Otherwise, $\mathfrak{q}$ is said to be of the {\it principal noetherian class} ".  \\

For example, denoting by $par_q$ the list of parametric jets up to order $q$, we have:  \\
\noindent
$ n=1 \Rightarrow n+1=2=2^1; \mathfrak{q}=({\chi}^2)=(\chi)^2 \Rightarrow par_1=\{y,y_x\} \Rightarrow dim_k(M)=2$.  \\
\noindent
$ n=2 \Rightarrow n+1=3<2^2=4; \mathfrak{q}=(({\chi}_2)^2, {\chi}_1{\chi}_2-({\chi}_1)^2)\subset ({\chi}_1,{\chi}_2)^2 \Rightarrow \\
par_2=\{y,y_1,y_2,y_{11}\} \Rightarrow dim_k(M)=2^2=4$. However, $\mathfrak{q'}=(({\chi}_2)^3, {\chi}_1{\chi}_2-({\chi}_1)^2) \Rightarrow par_3=\{ y,y_1,y_2,y_{11}, y_{22}, y_{111}\}  \Rightarrow dim_k(M')=6>4$.   \\
\noindent
$n=3 \Rightarrow n+1=4<2^3=8; \mathfrak{q}=(({\chi}_3)^2,{\chi}_2{\chi}_3-({\chi}_1)^2,({\chi}_2)^2) \subset ({\chi}_1,{\chi}_2,{\chi}_3)^2 \Rightarrow \\
par_3=\{y,y_1,y_2,y_3, y_{11},y_{12},y_{13},y_{111}\} \Rightarrow dim_k(M)=2^3=8$.  \\

Our aim is to explain this result in a modern language and to illustrate it by providing a similar example for $n=4$. The importance of such an example is that it allows for the first time to exhibit symbols which are $2,3,4$-acyclic without being involutive. Another interest of this example is that it has properties quite similar to the ones held by the system of conformal Killing equations which are still not known. For this reason, we have put all the examples at the end of the paper and each one is presented in 
a rather independent way though a few among them are quite tricky. \\

Meanwhile, the second purpose is to prove that the methods developped by Macaulay in order to study {\it unmixed polynomial ideals} are only particular examples of new formal differential geometric techniques that have been introduced recently in order to study {\it pure differential modules}. However these procedures are based on the formal theory of systems of ordinary differential (OD) or partial differential (PD) equations, in particular on a systematic use of the Spencer operator, and are still not acknowledged by the algebraic community.  \\
\hspace*{1cm}  \\

Commutative algebra, Localization, Macaulay inverse systems, Duality, \\ 
Spencer operator, Spencer cohomology, Involution.  \\

\noindent
{\bf 1) DIFFERENTIAL SYSTEMS   }\\

If $E$ is a vector bundle over the base manifold $X$ with projection $\pi$ and local coordinates $(x,y)=(x^i,y^k)$ projecting onto $x=(x^i)$ for $i=1,...,n$ and $k=1,...,m$, identifying a map with its graph, a (local) section $f:U\subset X \rightarrow E$ is such that $\pi\circ f =id$ on $U$ and we write $y^k=f^k(x)$ or simply $y=f(x)$. For any change of local coordinates $(x,y)\rightarrow (\bar{x}=\varphi(x),\bar{y}=A(x)y)$ on $E$, the change of section is $y=f(x)\rightarrow \bar{y}=\bar{f}(\bar{x})$ such that ${\bar{f}}^l(\varphi(x)\equiv A^l_k(x)f^k(x)$. The new vector bundle $E^*$ obtained by changing the {\it transition matrix} $A$ to its inverse $A^{-1}$ is called the {\it dual vector bundle} of $E$. Differentiating with respect to $x^i$ and using new coordinates $y^k_i$ in place of ${\partial}_if^k(x)$, we obtain ${\bar{y}}^l_r{\partial}_i{\varphi}^r(x)=A^l_k(x)y^k_i+{\partial}_iA^l_k(x)y^k$. Introducing a multi-index $\mu=({\mu}_1,...,{\mu}_n)$ with length $\mid \mu \mid={\mu}_1+...+{\mu}_n$ and prolonging the procedure up to order $q$, we may construct in this way, by patching coordinates, a vector bundle $J_q(E)$ over $X$, called the {\it jet bundle of order} $q$ with local coordinates $(x,y_q)=(x^i,y^k_{\mu})$ with $0\leq \mid\mu\mid \leq q$ and $y^k_0=y^k$. For a later use, we shall set $\mu+1_i=({\mu}_1,...,{\mu}_{i-1},{\mu}_i+1,{\mu}_{i+1},...,{\mu}_n)$ and define the operator $j_q:E \rightarrow J_q(E):f \rightarrow j_q(f)$ on sections by the local formula $j_q(f):(x)\rightarrow({\partial}_{\mu}f^k(x)\mid 0\leq \mid\mu\mid \leq q,k=1,...,m)$. Finally, a jet coordinate $y^k_{\mu}$ is said to be of {\it class} $i$ if ${\mu}_1=...={\mu}_{i-1}=0, {\mu}_i\neq 0$. \\

\noindent
{\bf DEFINITION 1.1}:  A {\it system} of PD equations of order $q$ on $E$ is a vector subbundle $R_q\subset J_q(E)$ locally defined by a constant rank system of linear equations for the jets of order $q$ of the form $ {\Phi}^{\tau} \equiv a^{\tau\mu}_k(x)y^k_{\mu}=0$. Its {\it first prolongation} $R_{q+1}\subset J_{q+1}(E)$ will be defined by the equations $ {\Phi}^{\tau} \equiv a^{\tau\mu}_k(x)y^k_{\mu}=0, d_i{\Phi}^{\tau} \equiv a^{\tau\mu}_k(x)y^k_{\mu+1_i}+{\partial}_ia^{\tau\mu}_k(x)y^k_{\mu}=0$ which may not provide a system of constant rank as can easily be seen for $xy_x-y=0 \Rightarrow xy_{xx}=0$ where the rank drops at $x=0$.\\

The next definition will be crucial for our purpose.\\

\noindent
{\bf DEFINITION 1.2}: A system $R_q$ is said to be {\it formally integrable} if the $R_{q+r}$ are vector bundles $\forall r\geq 0$ (regularity condition) and no new equation of order $q+r$ can be obtained by prolonging the given PD equations more than $r$ times, $\forall r\geq 0$.\\

Finding an intrinsic test has been achieved by D.C. Spencer in 1965 ([30]) along coordinate dependent lines sketched by M.Janet as early as in 1920 ([6]). The key ingredient, missing totally before the moderrn approach, is provided by the following definition.\\

\noindent
{\bf DEFINITION 1.3}: The family $g_{q+r}$ of vector spaces over $X$ defined by the purely linear equations $ a^{\tau\mu}_k(x)v^k_{\mu+\nu}=0$ for $ \mid\mu\mid= q, \mid\nu\mid =r $ is called the {\it symbol} at order $q+r$ and only depends on $g_q$.\\

The following procedure, {\it where one may have to change linearly the independent variables if necessary}, is the heart towards the next definition which is intrinsic even though it must be checked in a particular coordinate system called $\delta$-{\it regular} (See [14],[15] and [28],[29] for more details):\\

\noindent
$\bullet$ {\it Equations of class} $n$: Solve the maximum number ${\beta}^n_q$ of equations with respect to the jets of order $q$ and class $n$. Then call $(x^1,...,x^n)$ {\it multiplicative variables}.\\
$\bullet$ {\it Equations of class} $i$: Solve the maximum number of {\it remaining} equations with respect to the jets of order $q$ and class $i$. Then call $(x^1,...,x^i)$ {\it multiplicative variables} and $(x^{i+1},...,x^n)$ {\it non-multiplicative variables}.\\
$\bullet$ {\it Remaining equations equations of order} $\leq q-1$: Call $(x^1,...,x^n)$ {\it non-multiplicative variables}.\\

\noindent
{\bf DEFINITION 1.4}: A system of PD equations is said to be {\it involutive} if its first prolongation can be achieved by prolonging its equations only with respect to the corresponding multiplicative variables. The numbers ${\alpha}^i_q=m(q+n-i-1)!/((q-1)!(n-i)!)-{\beta}^i_q$ will be called {\it characters} and ${\alpha}^1_q\geq ... \geq {\alpha}^n_q $. For an involutive system, $(y^{{\beta}^n_q +1},...,y^m)$ can be given arbitrarily.  \\

Though the preceding description was known to Janet, he surprisingly never used it explicitly. In any case, such a definition is far from being intrinsic and the hard step will be achieved from the Spencer cohomology that will also play an important part in any explicit computation.\\

Let $T$ be the tangent vector bundle of vector fields on $X$, $T^*$ be the cotangent vector bundle of 1-forms on $X$ and ${\wedge}^sT^*$ be the vector bundle of s-forms on $X$ with usual bases $\{dx^I=dx^{i_1}\wedge ... \wedge dx^{i_s}\}$ where we have set $I=(i_1< ... <i_s)$. Also, let $S_qT^*$ be the vector bundle of symmetric q-covariant tensors. Finally, we may introduce the {\it exterior derivative} $d:{\wedge}^rT^*\rightarrow {\wedge}^{r+1}T^*:\omega={\omega}_Idx^I \rightarrow d\omega={\partial}_i{\omega}_Idx^i\wedge dx^I$ with $d^2=d\circ d\equiv 0$ in the {\it Poincar\'{e} sequence}:\\
\[  {\wedge}^0T^* \stackrel{d}{\longrightarrow} {\wedge}^1T^* \stackrel{d}{\longrightarrow} {\wedge}^2T^* \stackrel{d}{\longrightarrow} ... \stackrel{d}{\longrightarrow} {\wedge}^nT^* \longrightarrow 0  \]

In a purely algebraic setting, one has ([14],[28],[29],[30]):  \\

\noindent
{\bf PROPOSITION 1.5}: There exists a map $\delta:{\wedge}^sT^*\otimes S_{q+1}T^*\otimes E\rightarrow {\wedge}^{s+1}T^*\otimes S_qT^*\otimes E$ which restricts to $\delta:{\wedge}^sT^*\otimes g_{q+1}\rightarrow {\wedge}^{s+1}T^*\otimes g_q$ and ${\delta}^2=\delta\circ\delta=0$.\\

{\it Proof}: Let us introduce the family of s-forms $\omega=\{ {\omega}^k_{\mu}=v^k_{\mu,I}dx^I \}$ and set $(\delta\omega)^k_{\mu}=dx^i\wedge{\omega}^k_{\mu+1_i}$. We obtain at once $({\delta}^2\omega)^k_{\mu}=dx^i\wedge dx^j\wedge{\omega}^k_{\mu+1_i+1_j}=0$.\\
\hspace*{12cm} Q.E.D.  \\

The kernel of each $\delta$ in the first case is equal to the image of the preceding $\delta$ but this may no longer be true in the restricted case and we set:\\

\noindent
{\bf DEFINITION 1.6}: Let $B^s_{q+r}(g_q)\subseteq Z^s_{q+r}(g_q)$ and $H^s_{q+r}(g_q)=Z^s_{q+r}(g_q)/B^s_{q+r}(g_q)$ be respectively the coboundary space, cocycle space and cohomology space at ${\wedge}^sT^*\otimes g_{q+r}$ of the restricted $\delta$-sequence which only depend on $g_q$ and may not be vector bundles. The symbol $g_q$ is said to be s-{\it acyclic} if $H^1_{q+r}=...=H^s_{q+r}=0, \forall r\geq 0$, {\it involutive} if it is n-acyclic and {\it finite type} if $g_{q+r}=0$ becomes trivially involutive for r large enough. Finally, $S_qT^*\otimes E$ is involutive $\forall q\geq 0$ if we set $S_0T^*\otimes E=E$. \\

The following {\it Janet/Goldschmidt/Spencer criterion} is superseding the various procedures used by Macaulay, Janet or Gr\"obner ([6], [9], [28]):  \\

\noindent
{\bf THEOREM 1.7}: If $g_q$ is $2$-acyclic, $g_{q+1}$ is a vector bundle and ${\pi}^{q+1}_q:R_{q+1} \rightarrow R_q$ is an epimorphism, then $R_q$ is formally integrable. If moreover $g_q$ is involutive, then $g_q$ is {\it also} a vector bundle and $R_q$ is involutive.\\

\noindent
{\bf DEFINITION 1.8}: The {\it characteristic variety} $V$ of an involutive system of order $q$ is the algebraic set defined by the {\it radical} (care) of the polynomial ideal in $K[\chi]$ generated by the $m\times m$ minors of the {\it characteristic matrix} $(a^{\tau\mu}_k{\chi}_{\mu})$ where $\mid\mu\mid =q$. Of course, when $m=1$ we recover the radical of the symbol ideal of the ideal we started with and the involutive assumption is essential.\\

\noindent
{\bf 2) DIFFERENTIAL MODULES}  \\

We start this section with a rather unusual survey of {\it primary decomposition} in module theory, using a few technical results on {\it localization} which are not so well known ([1],[3],[8],[10],[16],[26]).  \\

Let $A$ be a commutative unitary integral domain, the best example being the ring $k[{\chi}_1,...,{\chi}_n]=k[\chi]$ of polynomials with coefficients in a field $k$. We denote as usual by $spec(A)$ the set of {\it proper} prime ideals of $A$ and by $max(A)$ the subset of maximal ideals. For any ideal $\mathfrak{a} \in A$, we introduce the family $Z(\mathfrak{a})=\{ \mathfrak{p}\in spec(A)\mid \mathfrak{p} \supset \mathfrak{a} \}$ of prime ideals containing $\mathfrak{a}$. We also introduce its {\it radical} by setting $rad(\mathfrak{a})=\{ a\in A \mid \exists r\in \mathbb{N}, a^r\in \mathfrak{a}\} \supseteq \mathfrak{a}$. We recall that $\mathfrak{q}$ is a {\it primary} ideal in $A$ if $ab\in \mathfrak{a}, a\notin \mathfrak{a} \Rightarrow b^r\in \mathfrak{a}$ for some integer $r\geq 1$ or simply $b\in rad(\mathfrak{a})$. In that case, $rad(\mathfrak{q})=\mathfrak{p} \in spec(A)$.  \\

\noindent
{\bf PROPOSITION 2.1}: If $0\rightarrow M' \stackrel{f}{\rightarrow } M \stackrel{g}{\rightarrow} M" \rightarrow 0$ is a short exact sequence, one has the following relations showing the importance of the radical: 
\[ ann(M) \subseteq ann(M')\cap ann(M") \Rightarrow rad(ann(M))=rad(ann(M'))\cap rad(ann(M")) \]
 
\noindent
{\bf DEFINITION 2.2}: A subset $S\subset A$ is said to be {\it multiplicatively closed} if $1\in S$ and $s,t\in S \Rightarrow st \in S$. If $M$ is a module over $A$, we shall set $t_S(M)=\{ x\in M \mid \exists s\in S, sx=0 \}$ and we have the exact sequence $0\rightarrow  t_S(M) \rightarrow M \stackrel{\theta}{\longrightarrow} S^{-1}M$ where $\theta = {\theta}_S=M \rightarrow S^{-1}M:x \rightarrow \frac{1}{1}x=\frac{1}{s} sx, \forall s\in S$ is the {\it localization} morphism over $S$. If $S=A-\mathfrak{p}$ with $\mathfrak{p} \in spec(A)$, we set as usual $S^{-1}M=M_{\mathfrak{p}}$. \\

\noindent
{\bf DEFINITION 2.3}: We say that $\mathfrak{p} \in spec(A)$ is an {\it associated prime} of $M$ if there exists $x\in M$ such that $ann(x)=\mathfrak{p}$. The set of all such associated primes will be denoted by $ass(M)$. \\

\noindent
{\bf LEMMA 2.4}: If $A$ is also noetherian and $M$ is a module over $A$, then $M \neq 0  \Leftrightarrow ass(M) \neq \emptyset $. Accordingly, the set of zero divisors for $M$ in $A$ is the union of all the associated primes and the set of non-zerodivisors is thus $S=A-{\cup}_{{\mathfrak{p}}_i \in ass(M)}{\mathfrak{p}}_i= {\cap}_{{\mathfrak{p}}_i \in ass(M)}(A-{\mathfrak{p}}_i)={\cap}_{{\mathfrak{p}}_i \in ass(M)} S_i$. \\

\noindent
{\bf DEFINITION 2.5}: If $P$ is a module over $A$, we say that $P$ is a {\it prime module} if we have $ax=0, x\neq 0 \Rightarrow a\in ann(P)$.  More generally, if $Q$ is a module over $A$, we say that $Q$ is a {\it primary module} if $ax=0, x\neq 0 \Rightarrow a\in rad(ann(Q))$. We shall say that a submodule $N\subset M$ is {\it prime } ({\it primary}) " {\it in} " $M$ if $M/N$ is a prime (primary) module. Accordingly, if $Q$ is a primary module, then $\mathfrak{q}=ann(Q)$ is a primary ideal in $A$ with $\mathfrak{p}=rad(\mathfrak{q})\in spec(A)$ and we shall say that $Q$ is $\mathfrak{p}$-primary. Similarly, if $P$ is a prime module with $\mathfrak{p}=ann(P)\in spec(A)$, we shall say that $P$ is $\mathfrak{p}$-prime.\\

\noindent
{\bf LEMMA 2.6}: Any submodule of a primary (prime) module is again a primary (prime) module. Moreover, if $Q_1$ and $Q_2$ are two $\mathfrak{p}$-primary modules, then $Q_1 \oplus Q_2$ is also a $\mathfrak{p}$-primary module.  \\

In actual practice, $Q$ is a primary ideal if and only if $\mathfrak{q} \subseteq ann(x) \subset \mathfrak{p}, \forall 0\neq x\in Q$ and $P$ is a prime module if and only if $\mathfrak{p}=ann(x), \forall 0\neq x\in P$. One can prove: \\

\noindent
{\bf PROPOSITION 2.7}: If $M$ is a finitely generated module over $A$, then $M$ is a primay module if and only if $ass(M)=\{ \mathfrak{p} \}$ consists of one element $\mathfrak{p} \in spec(A)$ only and $M$ is a prime module if and only if $ass(M)=\{ \mathfrak{p}\}$ consists of one element $\mathfrak{p}\in spec(A)$ only such that $ann(M)=\mathfrak{p}$. \\

A next step is provided by the following definition:ÊÊ\\

\noindent
{\bf DEFINITION 2.8}: We say that a submodule $N\subset M$ is {\it reducible} " {\it in} " $M$ if one can find two submodules $N_1\neq N$ and $N_2\neq N$ such that $N=N_1cap N_2$. Otherwise, $N$ is called {\it irreducible} in $M$. Passing to the residue, we have an embedding $0 \rightarrow M/N \stackrel{(p_1,p_2)}{\longrightarrow} (M/N_1) \oplus (M/N_2) $ where $p_1:M \rightarrow M/N_1$ and $p_2:M\rightarrow M/N_2$ are the canonical projections. Accordingly, we say that a module $M$ is a {\it reducible module} if one can find two modules $Q_1$ and $Q_2$ both with a monomorphisms $0 \rightarrow M \stackrel{(p_1,p_2)}{\longrightarrow} Q_1 \oplus Q_2$ where $p_1:M \rightarrow Q_1$ and $p_2:M \rightarrow Q_2$ are two epimorphisms. Otherwise, $M$ is called an {\it irreducible module}.  \\

If $M'\subset M$, defining $Q'_i=p_i(M')$ for $i=1,2$, we obtain by restriction: \\

\noindent
{\bf LEMMA 2.9}: Any submodule of a reducible module is again reducible.  \\

 \noindent
 {\bf PROPOSITION 2.10}: Any non-zero irreducible module $M$ over a noetherian ring $A$ is primay.  \\
 
 \noindent
 {\bf Proof}: Let $M$ be an irreducible module. As $M\neq 0 \Rightarrow ass(M)\neq \emptyset$, let us suppose that ${\mathfrak{p}}_1, {\mathfrak{p}}_2 \in spec(A)$ with ${\mathfrak{p}}_1 \neq {\mathfrak{p}}_2$ are such that ${\mathfrak{p}}_1 , {\mathfrak{p}}_2 \in ass(M)$. Then one could find two elements $x_i\in M$ with ${\mathfrak{p}}_i=ann(x_i)$ for $i=1,2$. Defining $Q_1$ and $Q_2$ by the long exact ker/coker sequences $0 \rightarrow {\mathfrak{p}}_i \rightarrow A \stackrel{x_i}{\rightarrow} M \stackrel{p_i}{\rightarrow} Q_i \rightarrow 0$ for $i=1,2$ with maps $x_i:a \rightarrow ax_i$, we obtain a map $M \stackrel{(p_1,p_2)}{\longrightarrow} Q_1 \oplus Q_2$. If $x\in ker((p_1,p_2))=ker(p_1) \cap ker(p_2)$, we have $x_i \in ker(p_i) \Rightarrow ann(x)={\mathfrak{p}}_i$ and this is impossible unless $x=0$ because the two associated primes are supposed to be different. Hence, the map $(p_1,p_2)$ is a monomorphism, a result showing that $M$ is reducible. It follows that $ass(M)=\mathfrak{p}$ and $M$ is a primary module. \\ 
 \hspace*{12cm}    Q.E.D.  \\
 
\noindent
{\bf DEFINITION 2.11}: We say that a module $M$ can be {\it nicely embedded} into a direct sum of modules if one can find modules $Q_i$ and epimorphisms $p_i:M \rightarrow Q_i$ for $i=1, ..., t$ such that we have a monomorphism $0 \rightarrow M \stackrel{(p_1,...,p_t)}{\longrightarrow} Q_1 \oplus ...  \oplus Q_t$.  \\

The following proposition essentially uses the two assumptions of the nice embedding, namely $p_1, ... , p_t$ are epimorphisms and $(p_1, ... , p_t)$ is a monomorphism (For more details, see [3], Chapter IV, \S 1, Exercise 11 and [16]).\\

\noindent
{\bf PROPOSITION 2.12}: If $M$ is a finitely generated module, there is an exact sequence $0 \rightarrow  M \stackrel{({\theta}_1, ... , {\theta}_t)}{\longrightarrow} {\oplus}_{{\mathfrak{p}}_i\in ass(M)} {M}_{{\mathfrak{p}}_i}$ and $M$ can be thus nicely embedded into a direct sum $Q_1 \oplus ... \oplus Q_t$ of modules. In this case, we have $ann(M)=ann(Q_1)\cap ... \cap ann(Q_t)$.  \\

\noindent
{\bf Proof}: With $ass(M)=\{ {\mathfrak{p}}_1, ... , {\mathfrak{p}}_t \}$ with ${\mathfrak{p}}_i\in spec(A)$ for $i=1, ... ,t$, let us introduce the canonical localization morphisms 
${\theta}_i: M \longrightarrow {M}_{{\mathfrak{p}}_i}=S^{-1}_iM$ with $S_i=A-{\mathfrak{p}}_i$. Setting $N_i=ker({\theta}_i)=t_{S_i}(M), Q_i=M/N_i\simeq im({\theta}_i)$, we can introduce the canonical projections $p_i$ in the short exact sequences $0 \rightarrow N_i \rightarrow M \stackrel{p_i}{\longrightarrow} Q_i \rightarrow 0$. Finally, introducing $N=N_1\cap ... \cap N_t=ker((p_1,...,p_t))\subset M$, we have only to prove that $N=0$. Otherwise, if $0\neq x \in N$, we may find $s_i\in S_i$ such that $s_ix=0$ for $i=1, ..., t$. However, according to the last lemma, $N\neq 0 \Rightarrow ass(N)\neq \empty$. Let thus $\mathfrak{p} \in ass(N)\subset ass(M) \Rightarrow \mathfrak{p}={\mathfrak{p}}_i$ for some $i$ and thus $\exists x\in N\subset N_i, {\mathfrak{p}}_i=ann(x)$ but $s_i\notin {\mathfrak{p}}_i$, a result leading to a contradiction.      \\
Finally, $M \subset {\oplus}^t_{i=1}Q_i \Rightarrow ann(M) \supseteq ann({\oplus}^t_{i=1}Q_i)={\cap}^t_{i=1}ann(Q_i)$. However, as $p_i$ is an epimorphism, we have 
$ann(M)\subseteq ann(Q_i) \Rightarrow ann(M) \subseteq {\cap}^t_{i=1}ann(Q_i)$ and thus $ann(M)=ann(Q_1)\cap ... \cap ann(Q_t)$.  \\
\hspace*{12cm}  Q.E.D.     \\

\noindent
{\bf THEOREM 2.13}: ({\it Existence of a primary decomposition}) Any submodule of a finitely generated module $M$ over $A$ admits a primary decomposition in $M$, that is can be expressed as the intersection of a finite number of submodules of $M$ that are primary in $M$ and correspond to different prime ideals $\{ {\mathfrak{p}}_1, ..., {\mathfrak{p}}_t \}$.  \\

Applying the theorem to the zero submodule of $M$, we may deal with $M$ alone and obtain the following generalization of the chinese remainder theorem:  \\

\noindent
{\bf COROLLARY 2.14}: ({\it Existence of a primary summation}) Any finitely generated module $M$ over $A$ can be nicely embedded into a direct sum of primary modules.\\

\noindent
{\bf THEOREM 2.15}: ({\it First uniqueness theorem}) $ass(M)=\{{\mathfrak{p}}_1,...,{\mathfrak{p}}_t \}$. \\

The following theorem explains the importance of {\it pure modules}, that is to say modules with unmixed annihilators, and a module is said to be $r$-{\it pure} when each isolated component is of codimension equal to $r$ when $A$ is a polynomial ring or a ring of ordinary or partial differential operators with constant coefficients.  \\

\noindent
{\bf THEOREM 2.16}: ({\it Second uniqueness theorem}) If a module $M$ can be nicely embedded into a reduced primary summation ${\oplus}^t_{i=1}Q_i$ where each 
$Q_i$ is a ${\mathfrak{p}}_i$-primary module, with $\{{\mathfrak{p}}_1,...,{\mathfrak{p}}_s \}$ the minimal primes of $ass(M)$, then {\it only} $Q_1, ..., Q_s$ are uniquely determined by $M$ through the images of the respective localization morphisms ${\theta}_1, ... ,{\theta}_s$. \\

{\it A basic idea in module theory is to look for the greatest semi-simple submodule of a given module}. For this, if $\mathfrak{m}\in max(A)\cap ass(M)$, then one can find a finite number of elements $x,y,...\in M$ killed by $\mathfrak{m}$. Accordingly, the map $x:A\rightarrow M:a\rightarrow ax$ has kernel $\mathfrak{m}$ and $A/\mathfrak{m}\simeq Ax\subseteq M$ is a simple module, like $Ay$ which may eventually be different and so on. \\

\noindent
{\bf DEFINITION  2.17}:  The direct sum $Ax\oplus Ay\oplus ...$ is called the {\it socle} of $M$ at $\mathfrak{m}$ and denoted by ${soc}_{\mathfrak{m}}(M)$. These simple components are called {\it isotypical} as they are all isomorphic to $A/\mathfrak{m}$. The {\it socle} of $M$ is $soc(M)=\oplus {soc}_{\mathfrak{m}}(M)$ for $\mathfrak{m}\in  max(A)\cap ass(M)$. It is the largest semi-simple submodule of $M$. \\

We notice that the double condition on the direct sum is essential as we need not only a submodule ($\mathfrak{m}\in ass(M)$) but {\it also} a simple module ($\mathfrak{m}\in max(A)\subseteq spec(A)$). Finally, $M$ is semi-simple if $M=soc(M)$ and $soc(M)=0$ if $M$ has no simple submodule, like the $\mathbb{Z}$-module $\mathbb{Z}$. Of course, when $M$ is a differential module, we obtain at once $soc(M)\subseteq t_{n-1}(M)\subseteq M$ (See Proposition 2.28). \\

\noindent
{\bf DEFINITION 2.18}: The {\it radical} of a module $R$ is the submodule $rad(R)$ which is the intersection of all the maximum {\it proper} submodules of $R$. If $rad(R)=0$, for example if $R$ is simple, we say that $R$ has no radical. If $R$ has no proper maximum submodule, then $rad(R)=R$.\\

\noindent
{\bf LEMMA 2.19}: $rad(R)$ is the intersection of all the kernels of the nonzero morphisms $R\rightarrow S$ where $S$ is a simple module.\\

\noindent
{\bf DEFINITION 2.20}: The {\it top} of the module $R$ is the semi-simple module defined by the short exact sequence $0\rightarrow rad(R)\rightarrow R \rightarrow top(R)\rightarrow 0$. It can also be defined as the largest quotient of $R$ that is a direct sum of simple modules.\\

\noindent
{\bf LEMMA 2.21}: If $R\neq 0$ is finitely generated, then $rad(R)\neq R$.\\

\noindent
{\bf PROPOSITION 2.22}: Tensoring by $R$ the short exact sequence $0 \rightarrow \mathfrak{m} \rightarrow A \rightarrow A/ \mathfrak{m} \rightarrow 0$, we get the short exact sequence $ 0 \rightarrow \mathfrak{m} R \rightarrow R \rightarrow (A/ \mathfrak{m})\otimes R  \rightarrow 0$ and thus the short exact sequence 
$ 0 \rightarrow \cap\mathfrak{m} R \rightarrow R \rightarrow top(R) \rightarrow 0$ coming from the {\it chinese remainder theorem} where the intersection is taken on all $\mathfrak{m} \in ass(M)\cap max(A)$.  \\

Let $K$ be a differential field with $n$ commuting derivations ${\partial}_i$ for $i=1,...,n$ and subfield of constants $k\subseteq K$. The (noncommutative) ring 
$D=K[d_1,...,d_n]=K[d]$ of differential operators with coefficients in $K$ has elements of the form $P=a^{\mu}d_{\mu}$ such that $\mid \mu \mid < \infty$ with an implicit summation on multi-indices and we have $d_ia=ad_i + {\partial}_ia, \forall a\in K$.Then $D$ is filtred by the order $q$ of operators and we have $K=D_0\subset D_1\subset ... \subset D_{\infty}=D$. Accordingly, $D$ is generated by $K$ and $T=D_1/D_0$ with $D_1=K \oplus T$ if we identify an element $\xi={\xi}^id_i\in T$ with the vector field $\xi={\xi}^i(x){\partial}_i$ of differential geometry, but with ${\xi}^i\in K$ now. If we introduce {\it differential indeterminates} $y=(y^1,...,y^m)$, we may extend $d_iy^k_{\mu}=y^k_{\mu+1_i}$ to ${\Phi}^{\tau}\equiv a^{\tau\mu}_ky^k_{\mu}\stackrel{d_i}{\longrightarrow} d_i{\Phi}^{\tau}\equiv a^{\tau\mu}_ky^k_{\mu+1_i}+{\partial}_ia^{\tau\mu}_ky^k_{\mu}$ for $\tau=1,...,p$. Therefore, setting $Dy^1+...+Dy^m=Dy\simeq D^m$ and calling $I=D\Phi\subset Dy$ the {\it differential module of equations}, we obtain by residue the {\it differential module} or $D$-{\it module} $M=Dy/D\Phi$, denoting the residue of $y^k_{\mu}$ by ${\bar{y}}^k_{\mu}$ when there can be a confusion. Introducing the two free differential modules $F_0\simeq D^{m_0}, F_1\simeq D^{m_1}$, we obtain equivalently the {\it free presentation} $F_1\stackrel{{\cal{D}}}{\longrightarrow} F_0 \rightarrow M \rightarrow 0$ of order $q$ when $m_0=m, m_1=p$ and ${\cal{D}}=\Phi \circ j_q$ {\it is acting by composition on the right}. Similarly, the module $M$ is filtred by the order $q$ of the linear combinations of the jet coordinates $y_q$ in $D_qy$ allowing to describe elements of $M$ by residue and we have the {\it inductive limit} $M_0\subseteq M_1\subseteq ... \subseteq M_q\subseteq ... \subseteq M_{\infty}=M$ with $d_iM_q\subseteq M_{q+1}$ but we notice that 
$D_rD_q=D_{q+r} \Rightarrow D_rM_q= M_{q+r}, \forall q,r\geq 0 \Rightarrow M=DM_q, \forall q\geq 0$ {\it in this particular case} according to the following commutative and exact diagram where $F\simeq D^m\Rightarrow F_q\simeq D^m_q$ is a free differential module and $I_q=I \cap F_q$: \\
\[   \begin{array}{rcccccl}
  & 0 & & 0 & & 0 & \\
   & \downarrow & & \downarrow & & \downarrow & \\
   0 \rightarrow & I_q & \rightarrow & F_q & \rightarrow & M_q & \rightarrow 0  \\
 & \downarrow & & \downarrow & & \downarrow & \\
 0 \rightarrow & I & \rightarrow & F & \rightarrow & M & \rightarrow 0
 \end{array}  \]   
It also follows from noetherian arguments and formal integrability or involution that $D_ rI_q=I_{q+r}, \forall r\geq 0$ though we have in general $D_rI_s\subseteq I_{r+s}, \forall r\geq 0, \forall s<q$ only. For example, when $k={\mathbb{Q}}$, the system $y_{33}=0, y_{13}-y_2=0$ is defining a module $M$ having the finite free presentation  $D^2 \rightarrow D \rightarrow M \rightarrow 0$. In order to determine $M_2$, one has to take the residue of $D_2$ with respect to the vector space $I_2=k y_{33} + k y_{23} + k y_{22} + k (y_{13}-y_2)$ which is the intersection of $D_2$ with the image $I$ of the presentation morphism $D^2 \rightarrow D: (P,Q) \rightarrow P y_{33}+Q(y_{13}-y_2)$, a result showing why it is so important to start with a formally integrable/involutive operator. We invite the reader to treat similarly the third order system $y_{111}=0, y_{22}=0, y_{12}=0$ which is neither formally integrable nor involutive with a strict inclusion $D_1I_2\subset I_3$.\\

\noindent
{\bf DEFINITION 2.23}: We define the {\it system} $R=hom_K(M,K)=M^*$ and set $R_q=hom_K(M_q,K)=M_q^*$ as the {\it system of order q} in order to have now the {\it projective limit} $R=R_{\infty}\rightarrow ... \rightarrow R_q \rightarrow ... \rightarrow R_1 \rightarrow R_0$. Taking into account the differential geometric framework of Spencer ([14-17],[30]), if a system of PD equations of order $q$ is given as before, then $f_q\in R_q:y^k_{\mu}\rightarrow f^k_{\mu}\in K$ with ${\Phi}^{\tau} \rightarrow a^{\tau\mu}_kf^k_{\mu}=0$ defines a {\it section at order} $q$ and we set $f_{\infty}=f\in R$ for a {\it section}. It is only when the field of constants $k$ is used that we can speak about a formal power series solution. The reader will recognize the so-called {\it inverse systems} of Macaulay ([9], \S 59, p 67) but in a quite different and much more general intrinsic framework.\\

\noindent
{\bf DEFINITION 2.24}: A {\it modular equation} $E\equiv f^k_{\mu}a^{\mu}_k=0$ of order $q$ with $0\leq k\leq m, 0\leq \mid\mu\mid\leq q$ is just a way to write down a section $f_q\in R_q$ by using an implicit summation with {\it formal coefficients}. Of course, infinite summations may also be considered, the simplest example being that of the exp™nential $e^x$ with $f=(1,1,1, ...)$. The procedure is absolutely similar to the case $m=1,K=k$ where one uses the purely formal power series {\it notation} $\sum f_{\mu}\frac{x^{\mu}}{\mu !}$ for writing down a {\it section}, even though the variable $x$ has absolutely no meaning in the module framework. Finally, as noticed by Macaulay, if one considers the set of modular equations at order $q$ as a homogeneous linear system for the unknowns $a^{\mu}_k$ at order $q$, then of course the given coefficients $a^{\tau\mu}_k$ form a basis of solutions linearly independant over $K$ and indexed by $\tau$. This is the reason for which we have chosen a similar notation.\\

For linear systems with constant coefficients and one unknown, Macaulay has found by induction the following striking theorem needing eventually a linear change of the polynomial indeterminates/coordinates ([9], \S 47, p 48 and \S 58, p 65). Of course not a word is left for more than one unknown, a situation where only the previous formal methods can be used. \\

\noindent
{\bf THEOREM 2.25}: The number of independent modular equations/parametric jet coordinates of strict order $q$ in the system provided by an ideal of the principal class $\mathfrak{a}=(P_1, ..., P_r)$ of rank $r$ with homogeneous polynomials $P_1,...,P_r$ of degree $l_1,...,l_r$, that is such that the corresponding differential module $M$ is $r$-pure with $cd(M)=r$, is equal to the coefficient of $x^q$ in the series $(1-x^{l_1})...(1-x^{l_r})(1-x)^{-n}$. In particular, if $r=n$ and $l_1=...=l_n=2$, we obtain the series $(1-x^2)^n(1-x)^{-n}=(1+x)^n$ showing that the total number of parametric jets up to order $n$ is $2^n$ indeed.  \\

\noindent
{\bf COROLLARY 2.26}: The same formula is still valid when $P_1,...,P_r$ are no longer homogeneous but their respective leading terms of highest degree $l_1,...,l_r$ are also generating an (homogeneous) ideal of the principal class of rank $r$ (See Example 3 and the twisted cubic of Example 6 for fine counterexamples but see Example 5 and in particular the third case of Example 6 for fine examples).  \\

The following nontrivial proposition generalizes the results of Macaulay to arbitrary systems with variable coefficients because $K$ is a $D$-module with the standard action $(D,K)\rightarrow K: (d_i,a)\rightarrow {\partial}_ia$. However, it is not evident, at first sight, to endow $M^*$ with a structure of left $D$-module in general, unless $D$ is a commutative ring, that is $K=cst(K)=k$ ([2], Theorem 1.3.1, 21, [18], Theorem 3.89, p 487).\\

\noindent
{\bf PROPOSITION 2.27}: When $M$ is a left $D$-module, then $R$ is a left $D$-module too.\\

\noindent
{\bf Proof}: As $D$ is generated by $K=D_0$ and $T=D_1/D_0$, let us define:\\
\[   (af)(m)=af(m)=f(am) \hspace{1cm} \forall a\in K, \forall m\in M\]
\[   (\xi f)(m)=\xi f(m)-f(\xi m)  \hspace{1cm}  \forall \xi =a^id_i\in T, \forall m\in M  \]
It is easy to check that $d_ia=ad_i+{\partial}_ia$ in the operator sense and that $\xi\eta -\eta\xi =[\xi,\eta]$ is the standard bracket of vector fields. We finally get $(d_if)^k_{\mu}=(d_if)(y^k_{\mu})={\partial}_if^k_{\mu}-f^k_{\mu +1_i}$ that is {\it exactly} the {\it Spencer operator} ([17],[18],[24]). We notice that the restriction of $d_i$ to $g_{q+1}$ is {\it minus} ${\delta}_i$, giving rise to the Spencer map $\delta = dx^i\otimes {\delta}_i:g_{q+1} \rightarrow T^*\otimes g_q$ with the corresponding Spencer $\delta$-cohomology. We have $d_id_j=d_jd_i=d_{ij}, \forall i,j=1,...,n$ because $(d_id_jf)^k_{\mu}={\partial}_{ij}f^k_{\mu}-{\partial}_if^k_{\mu +1_j}-{\partial}_jf^k_{\mu +1_i} +f^k_{\mu +{1_i}+{1_j}}$ and $d_iR_{q+1}\subseteq R_q$, a result leading to $d_iR\subset R$  and a well defined operator $R\rightarrow T^*{\otimes}_K R:f \rightarrow dx^i {\otimes} d_if$. 
Alternatively and in a coherent way with differential geometry, if we have a linear system ${\Phi}^{\tau}=0$ defining $R_q$ and its first prolongation ${\Phi}^{\tau}=0,d_i{\Phi}^{\tau}=0 $ defining $R_{q+1}$, a section $f_{q+1}\in R_{q+1}$ over $f_q\in R_q$ satisfies both $a^{\tau\mu}_kf^k_{\mu}=0$ and $a^{\tau\mu}_kf^k_{\mu +1_i}+{\partial}_ia^{\tau\mu}_kf^k_{\mu}=0$ as equalities in $K$ with $0\leq \mid\mu\mid \leq q$. Applying ${\partial}_i$ to the first and substracting the second, we get $a^{\tau\mu}_k({\partial}_if^k_{\mu}-f^k_{\mu +1_i})=0$. Accordingly, we obtain $ f_{q+1}\in R_{q+1}\stackrel{d_i}{\longrightarrow}d_if_{q+1}\in R_q$ and thus: \\ 
\[  E\equiv f^k_{\nu}a^{\nu}_k=0  \stackrel{d_i}{\longrightarrow} d_iE\equiv ({\partial}_if^k_{\mu}-f^k_{\mu +1_i})a^{\mu}_k=0, \forall f\in R   \]
but $d_iE$ is of order $q$ with $0\leq\mid\mu\mid\leq q$ whenever $E$ is of order $q+1$ with $0\leq\mid\nu\mid\leq q+1$. When $K=k$, the partial derivative disappears and we recognize, {\it exactly but up to sign},  the operator of Macaulay ([9], \S 60, p 69, [18]). For this reason and unless mentioned explicitly, {\it in this specific situation only}, we shall change the sign of the Spencer operator in order to agree with Macaulay.\\
\hspace*{12cm}   Q.E.D.  \\

We now provide a few results originating from the study of differential modules in {\it algebraic analysis} (See [2],[7],[11],[12],[13],[17],[27] for more details). \\

If $m\in M$, then the differential submodule $Dm\subset M$ is defined by a system of OD/PD equations for one unknown only and we may look for its codimension $cd(Dm)$. In the commutative case, looking at the annihilators, we get $ann(M)\subset ann(Dm)$. In particular, if $M$ is primary its annihilator is a primary ideal $\mathfrak{q}$ with radical $\mathfrak{p}$ and we have $\mathfrak{q} \subseteq ann(Dm) \subseteq \mathfrak{p}, \forall m\in M$ as a possible characterisation. Accordingly, if $M$ is prime, then $ann(Dm)=\mathfrak{p}, \forall m\in M$.\\

\noindent
{\bf PROPOSITION 2.28}: $t_r(M)=\{m\in M\mid cd(Dm)>r\}$ is the greatest differential submodule of $M$ having codimension $>r$. It does not depend on the presentation and thus on the filtration of the module $M$. We have the nested chain of $n$ differential submodules:  
\[  0 =t_n(M) \subseteq t_{n-1}(M) \subseteq ... \subseteq t_1(M) \subseteq t_0(M)=t(M) \subseteq M   \]

Thanks to its implementation in [23-25], this proposition is essential for the use of computer algebra. \\

\noindent
{\bf PROPOSITION 2.29}: $cd(M)=cd(V)=r \Leftrightarrow {\alpha}^{n-r}_1\neq 0, {\alpha}^{n-r+1}_1= ... ={\alpha}^n_1=0 \Leftrightarrow t_r(M)\neq M,t_{r-1}(M)= ... =t_0(M)=t(M)=M$. \\

We may therefore define ([2],[7],[16]):\\

\noindent
{\bf DEFINITION 2.30}: $M$ is r-{\it pure} $\Leftrightarrow t_r(M)=0, t_{r-1}(M)=M, \Leftrightarrow cd(Dm)=r, \forall m\in M$. In particular, $M$ is 0-pure iff $t(M)=0$. Otherwise, if $cd(M)=r$ but $M$ is not pure, we shall call $M/t_r(M)$ the {\it pure part} of $M$.\\

The following result using a kind of {\it relative localization} generalizes for an arbitrary $m$ the similar ones first obtained by Macaulay for $m=1$ ([9], \S 82) and provides a technical test linking purity and involution, both with an effective construction of the two previous Propositions. From now on we shall only consider the constant coefficient situation, considering ${\chi}_1, ..., {\chi}_{n-r}$ just like parameters while setting $k'=k({\chi}_1,...,{\chi}_n)$ ([9], \S 77, p 86).\\

\noindent
{\bf THEOREM 2.31}: If $cd_D(M)=r$ one has the exact sequence:\\
\[    0 \longrightarrow t_r(M) \longrightarrow M \longrightarrow k({\chi}_1,...,{\chi}_{n-r})\otimes M  \]

\noindent
{\bf COROLLARY 2.32}: ([9], \S 41) $M$ is r-pure $\Leftrightarrow  0\rightarrow M \rightarrow k({\chi}_1,...,{\chi}_{n-r})\otimes M$ is exact.\\

\noindent
{\bf REMARK 2.33}: In actual practice it is important to notice that the relative localization {\it kills} the PD equations of class 1 up to class $(n-r-1)$ ({\it care}) because of the compatibility conditions provided by the involutive assumption and that the equations of strict class $(n-r)$ ({\it care again}) provide the smallest nonzero character. Moreover, any prime or primary module is pure.\\

\noindent
{\bf REMARK 2.34}: If $cd(M)=r$ and $K=k$, then ${\alpha}^{n-r+1}_1=0,..., {\alpha}^n_1=0$ and a relative localization  brings the system to a {\it finite type} (zero symbol at high order) system in $(d_{n-r+1},...,d_n)$ over the field $k({\chi}_1,...,{\chi}_{n-r})$. Accordingly, there is a finite number of linearly independent  sections of the localized system and thus an equal finite number of modular equations as. We have also $ann_D(M)=ann_D(R)$. Indeed, as a representative of any element of $M$ can be written as a {\it finite} linear combination of parametric jets with coefficients in $k$, we have $M\subseteq M^{**}$ and thus $ann(M)\subseteq ann(M^*)\subseteq ann(M^{**})\subseteq ann(M) \Rightarrow ann(M)=ann(M^*)$. This result, {\it not true at all in the noncommutative situation}, generalizes the one of Macaulay ([9], \S 61, p 70) obtained when $m=1$. \\

\noindent
{\bf REMARK 2.35}: When $K=k$, then $R=M^*=hom_k(M,k)=hom_D(M,D^*)$ where $D^*=hom_k(D,k)$ is an injective module ([4], Proposition 11, p 18). There is a canonical differential structure on $R$ induced by the Spencer operator, namely $(d_if)^k_{\mu}={\partial}_if^k_{\mu} - f^k_{\mu + 1_i}=- f^k_{\mu + 1_i}$, which is {\it exactly} (up to sign) the one introduced by Macaulay ([9], \S 60, p 69). The crucial idea of Macaulay has been to notice that $top(R) $ is the dual of $soc(M)$ by means of this {\it duality}, a result allowing to use Nakayama's lemma (Krull-Azumaya theorem) for finding generating sections (formal solutions) of $R$, {\it on the crucial condition to have a finitely generated differential module} $R$ {\it over} $D$, the true reason for which Macaulay had to introduce the relative localization in such a way that $M$ and thus $R$ become finite dimensional differential vector spaces over $k$ ([18],[19]). Moreover, if $S$ is a simple submodule of $M$, we have $\mathfrak{m}=ann(S)\subseteq ann(S^*)\subseteq ann(S^{**})=ann(S) \in max(A)$ in the purely algebraic framework. Acordingly, the dual $S^*$ of a simple module $S$ isomorphic to $A/ \mathfrak{a}$ is again an isotypical simple module, because else it would have a proper factor module, the dual of which would be a proper submodule of $A / \mathfrak{m}$. \\

\noindent
{\bf 3) EXAMPLES}  \\

\noindent
{\bf EXAMPLE 1}: ([5], p526) If $A=\mathbb{Q}[{\chi}_1,{\chi}_2]$ and $M=A/\mathfrak{a}$ with a polynomial ideal $\mathfrak{a}=(({\chi}_1)^3,({\chi}_2)^2,{\chi}_1{\chi}_2)$, we have $rad(\mathfrak{a})=({\chi}_1,{\chi}_2)=\mathfrak{m} \in max (A)$ and $\mathfrak{a}$ is a primary ideal. The corresponding differential module $M$ over $D=k[d_1,d_2]$ is defined by the third order homogeneous involutive system $y_{222}=0, y_{122}=0, y_{112}=0, y_{111}=0, y_{22}=0, y_{12}=0$ but {\it both} ${\bar{y}}_2$ and ${\bar{y}}_{11}$ are killed by $\mathfrak{m}=(d_1,d_2)$. It follows that $soc(M)=soc_{\mathfrak{m}}(M)=D{\bar{y}}_2 \oplus D {\bar{y}}_{11}$ has two isotypical components isomorphic to $D/\mathfrak{m}\simeq k$. We have $M\simeq k\bar{y} + k {\bar{y}}_{1} + k{\bar{y}}_{2} + k{\bar{y}}_{11}$ with the action of $D$ described by $d_i{\bar{y}}_{\mu}={\bar{y}}_{\mu + 1_i}$. In the present case, we have $R=ka^0 + k a^1 + k a^2 + k a^{11} \Rightarrow \mathfrak{m} R = k a^0 + k a^1$ along with Proposition 2.22 and $R$ is thus generated by the {\it two} modular equations $ E_1\equiv a^2=0$ {\it and} $E_2 \equiv a^{11}=0$ according to Nakayama's lemma, as can be checked directly . \\

\noindent
{\bf EXAMPLE 2}: With $n=3$, let us consider the mixed homogeneous polynomial ideal: \\
\[\mathfrak{a}=(({\chi}_3)^2, {\chi}_2{\chi}_3, {\chi}_1{\chi}_3, {\chi}_1{\chi}_2)=({\chi}_3,{\chi}_1)\cap ({\chi}_3,{\chi}_2) \cap ({\chi}_1, {\chi}_2, {\chi}_3)^2={\mathfrak{p}}_1 \cap {\mathfrak{p}}_2 \cap (\mathfrak{m})^2\] 
as in ([18], Ex 4.10). We have the exact sequence $0 \rightarrow M \stackrel{({\theta}_1, {\theta}_2, {\theta}_3)}{\longrightarrow} M_{{\mathfrak{p}}_1}\oplus M_{{\mathfrak{p}}_2} \oplus M_{\mathfrak{m}}$ and set $N_i=ker ({\theta}_i)$ for $i=1,2$. Let us study the {\it specialization} defined by the short exact sequence $0 \rightarrow N_1 \rightarrow M \stackrel{p_1}{\longrightarrow} Q_1 \rightarrow 0 $ implying that $Q_1$ is defined by {\it more} PD equations than $M$ in the differential framework. As $k=\mathbb{Q}$, the system $R=M^*=hom_k(M,k)$ is defined by the second order PD equations $y_{33}=0, y_{23}=0, y_{13}=0, y_{12}=0$ but is not involutive and does not therefore provide any information on $R$ and thus on $M$. Changing $d_1$ to $ d_1-d_2$ while keeping $d_2$ and $d_3$ unchanged, we get the involutive system: \\
\[ \left\{   \begin{array}{l}
y_{33}=0 \\
y_{23}=0, \hspace{3mm} y_{22}-y_{12} =0 \\
y_{13}=0
\end{array}
\right. \fbox{ $ \begin{array}{lll}
1 & 2 & 3 \\
1 & 2 & \bullet \\
1 & \bullet & \bullet
\end{array} $ }      \]
with full classes $3$ and $2$ providing $cd(M)=2$ and thus $rk(\mathfrak{a})=2$. It is important to notice that we have now only one equation of class $1$ left instead of two in the original system of coordinates and this fact makes all the difference. However, localizing with respect to $k'=k({\chi}_1)$ gives ${\chi}_1 y_3=0$ and thus ${\bar{y}}_3$ is a torsion element (in fact the one killed by $\mathfrak{m}$), a result showing that $M$ is {\it not} $2$-pure. We invite the reader to transform this system into a first order system like in the next example in order to exhibit the corresponding torsion part. Finally, a polynomial in $S=A-\mathfrak{m}$ must have a non-vanishing term of degree zero and cannot therefore kill any element in $M$ which is defined by homogeneous polynomials. It follows that $t_S(M)=0$ and $Q=M/t_S(M)=M$ is surely not primary. On the contrary, identifying $k[\chi]$ and $k[d_1,d_2,d_3]$, we have $S_1=A-{\mathfrak{p}}_1=k+kd_2+kd_{22}+ ...$ and $N_1$ is generated by ${\bar{y}}_3, {\bar{y}}_1$. Accordingly, $Q_1$ is defined by the system $y_3=0, y_1=0$, that is $Q_1=P_1$ is a prime and thus ${\mathfrak{p}}_1$-primary module. Similarly, $Q_2=P_2$ is defined by the system $y_3=0, y_2=0$ and is thus also a prime and thus ${\mathfrak{p}}_2$-primary module. \\

\noindent
{\bf EXAMPLE 3}: With $n=3$, following Macaulay as in ([9], \S 38, p 40) where one can find the first intuition of formal integrability, let us first prove that the polynomial ideal $\mathfrak{q}=(({\chi}_1)^2, {\chi}_1{\chi}_3 - {\chi}_2)\in k[\chi]=A$ is primary while its radical is the prime ideal $\mathfrak{p}=({\chi}_1,{\chi}_2)\in spec(A)$ which is such that ${\mathfrak{p}}^2 \subset \mathfrak{q} \subset \mathfrak{p}$ ([16], Ex 1.188, p 127). Indeed, if $P,Q\in A$ with $PQ\in \mathfrak{q}\subset \mathfrak{p}$ and $Q\notin \mathfrak{p}\Rightarrow P\in \mathfrak{p}$, we just need to prove that $P\in \mathfrak{q}$. Passing modulo $\mathfrak{q}$, we may write $\bar{Q}=a({\chi}_3)+b({\chi}_3)\bar{{\chi}_1}$ with $a({\chi}_3)\neq 0$ because $ \bar{{\chi}_2}=\bar{{\chi}_1}{\chi}_3 $ as ${\chi}_3$ is arbitrary. Similarly, we have $\bar{P}=c({\chi}_3)\bar{{\chi}_1}$ and obtain therefore: \\
\[   \bar{P}\bar{Q}= c({\chi}_3)\bar{{\chi}_1}(a({\chi}_3)+b({\chi}_3)\bar{{\chi}_1})=a({\chi}_3)c({\chi}_3)\bar{{\chi}_1}=0       \]
We finally get $c({\chi}_3)=0$ because $a({\chi}_3)\neq 0$ and thus $\bar{P}=0 \Rightarrow P\in \mathfrak{q}$.\\
The corresponding second order system $R_2$ of PD equations is $y_{11}=0, y_{13}-y_2=0$ and is neither formally integrable nor involutive. In any case, $R^{(1)}_2$ is easily seen to be defined by $y_{11}=0, y_{12}=0, y_{13}-y_2=0$ because $d_3y_{11}-d_1(y_{13}-y_2)=y_{12}$ while $R^{(2)}_2$ is defined by $y_{11}=0, y_{12}=0, y_{22}=0, y_{13}-y_2=0$ because $d_{33}y_{11}-d_{13}(y_{13}-y_2)-d_2(y_{13}-y_2)=y_{22}$. In order to check that this last system is formally integrable and even involutive, let us effect the permutation $(1,2,3) \rightarrow (3,2,1)$ in order to get the new system $y_{33}=0, y_{23}=0, y_{22}=0, y_{13}-y_2=0$. The symbol is in $\delta$-regular coordinates with characters ${\alpha}^3_2=0, {\alpha}^2_2=0, {\alpha}^1_2=\alpha=3-1=2$ and thus full classes $3$ and $2$ providing $rank(\mathfrak{q})=2$ in the language of Macaulay or $cd(M)=2$ in the language of differential modules. Also $\mathfrak{q}$ is of {\it principal class} and thus {\it unmixed} as it can be presented by means of the basis $\{P_1\equiv {({\chi}_3})^2 , P_2={\chi}_1{\chi}_3-{\chi}_2\}$ with only $2$ generators satisfying the regularity condition of ([9], \S 48, p 49) because $(({\chi}_3)^2)$ is of rank $1$. We notice that the homogeneous symbol part, namely $y_{33}=0, y_{23}=0, y_{22}=0, y_{13}=0$ is of course of rank $2$ but that the homogeneous symbol part of the basis considered is $y_{33}=0, y_{13}=0$ which is only of rank $1$ in a coherent way with ([M], \S 49, p 50), a result proving that $cd(M)=2$. Contrary to the previous example, let us prove that $M$ is $2$-pure. For this, transforming the system into a first order system while using the same $\delta$-regular coordinates, we may set $z^1=y, z^2=y_1, z^3=y_2, z^4=y_3$ in order to obtain the first order involutive system:  \\
\[ \left\{   \begin{array}{l}
z^1_3- z^4=0, \hspace{3mm} z^2_3-z^3=0,\hspace{3mm} z^3_3=0,\hspace{3mm}  z^4_3=0 \\
z^1_2-z^3=0, \hspace{3mm} z^2_2-z^3_1=0, \hspace{3mm}  z^3_2=0, \hspace{3mm} z^4_2=0  \\
z^1_1-z^2=0, \hspace{3mm} z^4_1-z^3=0
\end{array}
\right. \fbox{ $ \begin{array}{lll}
1 & 2 & 3 \\
1 & 2 & \bullet \\
1 & \bullet & \bullet
\end{array} $ }      \]
with characters ${\alpha}^3_1=0, {\alpha}^2_1=0, {\alpha}^1_1=\alpha=4-2=2$ and thus full classes $3$ and $2$ again. Localizing with respect to $k'=k({\chi}_1)$, the localized PD equations of class $1$ provide the parametrization $z^2={\chi}_1z^1, z^3={\chi}_1z^4$. The corresponding differential module over $k[d_1]$ is thus torsion-free and $M$ is therefore $2$-pure. We have thus illustrated on this example the way to extend to systems of PD equations with many unknowns the methods of Macaulay which were only valid for systems of PD equations with one unknown. Finally, we have $\mathfrak{q} \subset \mathfrak{p}$ and both ideals have rank $2$ while $\mathfrak{q}$ is of the principal ideal class {\it and } primary, a result showing that $\mathfrak{p}$ is of the principal noetherian class ([9], \S 68, p 79). \\

\noindent
{\bf EXAMPLE 4}: With $n=3$, let us consider the homogeneous polynomial ideal $\mathfrak{a}=(({\chi}_3)^2, {\chi}_2{\chi}_3 - {\chi}_1{\chi}_3, ({\chi}_2)^2 - {\chi}_1{\chi}_2)\in k[\chi]$. This ideal is {\it not} primary because ${\chi}_2({\chi}_2 - {\chi}_1)\in \mathfrak{a}$ and each factor cannot belong to the ideal which is homogeneous. As for the other factor, it cannot belong to $rad(\mathfrak{a})=({\chi}_3, {\chi}_2({\chi}_2 - {\chi}_1))$ for the same reason. Introducing the primary ideals ${\mathfrak{q}}_1= (({\chi}_3)^2,{\chi}_2 - {\chi}_1)$ and ${\mathfrak{q}}_2=({\chi}_3,{\chi}_2)$, we obtain at once $\mathfrak{a}={\mathfrak{q}}_1\cap {\mathfrak{q}}_2$. \\
The corresponding differential module $M$ is defined by the homogeneous second order involutive system $y_{33}=0, y_{23} - y_{13}=0, y_{22} - y_{12}=0$ with full classes $3$ and $2$, a fact providing $cd(M)=2$ and giving a rank equal to $r=2$ for the ideal. The first non-zero character is $\alpha = {\alpha}^{n-r}_2={\alpha}^1_2
= 3$. The ideal $\mathfrak{a}$ is unmixed though not of principal class while $M$ is $2$-pure. With a slight abuse of language, the localized module of equations is the intersection of the localized primary ideal ${\mathfrak{q}}_1=(d_{33},d_2 - {\chi}_1)$ allowing to define the system $R_1$, with ${\mathfrak{q}}_2=(d_3,d_2)$ allowing to define the system $R_2$, in such a way that $R=R_1 \oplus R_2$ because $d_2 - (d_2 - {\chi}_1)={\chi}_1$ and thus ${\mathfrak{q}}_1 + {\mathfrak{q}}_2=k({\chi}_1)[d_2,d_3]=k'[d_2,d_3]$. Accordingly, the localized system $R$ defined by $y_{33}=0, y_{23} - {\chi}_1 y_3=0, y_{22} - {\chi}_1 y_2=0$ with respect to $k({\chi}_1)[d_2,d_3]$ is such that $par_2=\{ y, y_2,y_3 \}$ and thus $dim_{k'}(R)=3=\alpha$. On the other side, we have $dim_{k'} (R_1)=2$ while $dim_{k'}(R_2)=1$ and $dim_k(R)=dim_k(R_1) + dim_k(R_2)$ in a coherent way with the previous direct sum. \\

\noindent
{\bf EXAMPLE 5}: Following ([9], end of \S 71, p 81) with $n=3$, let us consider the primary polynomial ideal ${\mathfrak{m}}^2$ with $\mathfrak{m}=({\chi}_1,{\chi}_2,{\chi}_3)\in max(k[\chi])$ and the primary polynomial ideals $\mathfrak{q}=(({\chi}_3)^2 - ({\chi}_1)^2, {\chi}_2{\chi}_3 , ({\chi}_2)^2 -({\chi}_1)^2, {\chi}_1{\chi}_3, {\chi}_1{\chi}_2), {\mathfrak{q}}'=( ({\chi}_3)^2 - ({\chi}_1)^2, {\chi}_2{\chi}_3, ({\chi}_2)^2- ({\chi}_1)^2), {\mathfrak{q}}"=( ({\chi}_3)^2 - ({\chi}_1)^2, {\chi}_2({\chi}_3)^2, ({\chi}_2)^2- ({\chi}_1)^2)$. We have ${\mathfrak{q}}" \subset {\mathfrak{q}}' \subset \mathfrak{q}\subset {\mathfrak{m}}^2$ and all these primary ideals have the same radical $\mathfrak{m}$ but only ${\mathfrak{q}}' $ and ${\mathfrak{q}}"$ are ideals of the principal class in regular form. The corresponding second and third order systems are $R$ defined by $y_{33} - y_{11}=0, y_{23}=0, y_{22} - y_{11}=0, y_{13}=0, y_{12}=0$, $R'$ defined by $y_{33} - y_{11}=0, y_{23}=0, y_{22} - y_{11}=0$ and $R"$ defined by $y_{33} - y_{11}=0, y_{233}=0, y_{22} - y_{11}=0$ which are finite type with $g_3=0, g'_4=0,g"_5=0$ and $dim_k(g_2)=1, dim_k(g'_3)=1, dim_k(g"_4)=1$, a result leading to $dim_k(R)=1+3+1=5$, $dim_k(R')=1+3+3+1=8=2^3$ because the Hilbert function is $(1+x)^3=1+3x+3x^2+x^3$ as we already saw and $dim_k(R")=1+3+4+3+1=12$ because the Hilbert function is now $(1+x+x^2)(1+x)^2=1+3x+4x^2+3x^3+x^4$. Also, the reader may easily check that $R$ is generated by the single modular equation $E\equiv a^{11} + a^{22} + a^{33}=0$ because $R= ka^0 + ka^1 + ka^2 + ka^3 + k(a^{11} + a^{22} + a^{33}) \Rightarrow \mathfrak{m} R= k a^0 +ka^1 + ka^2 +ka^3$, while $R'$ is generated by the single modular equation $E'\equiv a^{111} + a^{122} + a^{133}=0$ providing for example $d_1E'\equiv a^{11} + a^{22} + a^{33}=0$ by applying the Spencer operator $d_1$ to $E'$ (up to sign). Similarly, $R"$ is generated by the single modular equation $E"\equiv a^{1113}+a^{1223}+a^{1333}=0$. Finally, $\mathfrak{q'}$ and ${\mathfrak{q}}"$ are defining ideals of the principal class in regular form according to Macaulay and $dim_k(R)<dim_k(R')<dim_k(R")$ though we have $rad(\mathfrak{q})=rad({\mathfrak{q}}')=rad({\mathfrak{q}}")=\mathfrak{m}$ and $\mathfrak{q}$ is not of the principal noetherian class according to Macaulay ([9], \S 68, p 78).\\

\noindent
{\bf EXAMPLE 6}: Following the study of curves defined by complete intersections as in ([8], p 137-139), we use the differential module approach on a few explicit examples defined in cartesian coordinates.\\

\noindent
$\bullet$  The twisted cubic $({\chi}_1=u^3, {\chi}_2=u^2, {\chi}_3=u)$ is defined by the polynomial ideal $({\chi}_3)^3-{\chi}_1=0,({\chi}_3)^2-{\chi}_2=0$. The corresponding third order system $R_3$ defined by $y_{333}-y_1=0, y_{33}-y_2=0$ is not formally integrable and cannot give any information on the underlying differential module $M$. Using crossed derivatives, we get the second order system $R^{(2)}_2$ defined by $y_{33}-y_2=0, y_{23}-y_1=0$ with similar comments. Finally, the second order system $R^{(3)}_2$ defined by $y_{33}-y_2=0, y_{23}-y_1=0, y_{22}-y_{13}=0$ is involutive with $\delta$-regular coordinates. The differential module $M$ is $2$-pure because the classes $2$ and $3$ are full. We have at once  $par_2=\{ y, y_1,y_2,y_3,y_{11}, y_{12}, y_{13}\}$ but the system is not homogeneous and we cannot compute the Hilbert function as Macaulay did. It is quite important to notice that the original system is defined by ($1$ equation of order 
$3$ + $1$ equation of order $2$) while the " {\it useful} " involutive final system (which is never quoted !) is defined by ($3$ equations of order $2$). As the characters are ${\alpha}^3_2=0, {\alpha}^2_2=0, {\alpha}^1_2=3$, we get $dim(g_0)=1, dim(g_q)=3, \forall q\geq 1$.\\

\noindent
$\bullet$  The curve ${\chi}_1=u^3, {\chi}_2=u^4, {\chi}_3=u^5$ must be defined by $3$ polynomials and does not provide therefore a complete intersection. Such an example is providing a $2$-pure differential module $M$ and has been fully treated in ([19], Example 5.12).  \\

\noindent
$\bullet$  The curve ${\chi}_1=u^6, {\chi}_2=u^5, {\chi}_3=u^4$ is providing a complete intersection defined by the two polynomials $({\chi}_3)^3-({\chi}_1)^2=0, ({\chi}_2)^2-{\chi}_1{\chi}_3=0$. The corresponding third order system $y_{333}-y_{11}=0, y_{22}-y_{13}=0$ is not formally integrable and the same comment as before may be done. Using one prolongation, we get the following $1+2+1=4$ equations of the third order where the classes have been separated:
\[ \left\{   \begin{array}{l}
y_{333}-y_{11}=0,  \\
y_{223}-y_{133}=0, y_{222}-y_{123}=0,  \\
y_{122}-y_{113}=0
\end{array}
\right. \fbox{ $ \begin{array}{lll}
1 & 2 & 3 \\
1 & 2 & \bullet \\
1 & \bullet & \bullet
\end{array} $ }      \]
The system $R_3$ thus obtained by adding the only second order equation provided is not involutive. Using one prolongation, we get the following $1+4+4=9$ equations of the fourth order where the classes have been separated:
\[ \left\{   \begin{array}{l}
y_{3333}-y_{113}=0,  \\
y_{2333}-y_{112}=0, y_{2233}-y_{111}=0, y_{2223}-y_{1233}=0, y_{2222}-y_{1133}=0, \\
y_{1333}-y_{111}=0, y_{1223}-y_{1133}=0, y_{1222}-y_{1123}=0, y_{1122}-y_{1113}=0
\end{array}
\right. \fbox{ $ \begin{array}{lll}
1 & 2 & 3 \\
1 & 2 & \bullet \\
1 & \bullet & \bullet
\end{array} $ }      \]
The system $R_4$ thus obtained by adding the second and third order equations already provided is involutive with full classes $3$ and $2$, providing therefore a $2$-pure differential module $M$ because localizing with respect to ${\chi}_1$ does not provide new equations/corresponding torsion elements ([9], \S 41, p 43). The characters are ${\alpha}^3_4=0, {\alpha}^2_4=0, {\alpha}^1_4=10-4=6$. Using one prolongation, we get the following $1+5+9=15$ equations of the fifth order where the classes have been separated:
\[ \left\{   \begin{array}{l}
y_{33333}-y_{1133}=0  \\
y_{23333}-y_{1123}=0, y_{22333}-y_{1113}=0, y_{22233}-y_{1112}=0, \\ 
y_{22223}-y_{1111}=0, y_{22222}-y_{11233}=0  \\
y_{13333}-y_{1113}=0, y_{12333}-y_{1112}=0, y_{12233}-y_{1111}=0,   \\
y_{12223}-y_{11233}=0, y_{12222}-y_{11133}=0, y_{11333}-y_{1111}=0, \\
y_{11223}-y_{11133}=0, y_{11222}-y_{11123}=0, y_{11122}-y_{11113}=0
\end{array}
\right. \fbox{ $ \begin{array}{lll}
1 & 2 & 3 \\
1 & 2 & \bullet \\
1 & 2 & \bullet \\
1 & \bullet & \bullet  \\
1 & \bullet & \bullet  \\
1 & \bullet & \bullet
\end{array} $ }      \]
The symbol at order $5$ is thus again involutive and the corresponding system $R_5$ finally obtained by adding the second, third and fourth order equations already provided is again involutive with full classes $3$ and $2$. Collecting all the results so far obtained, we obtain the list:  \\
 $  par_5=\{ y, y_1, y_2, y_3, y_{11}, y_{12}, y_{13}, y_{23}, y_{33}, y_{111}, y_{112}, y_{113}, y_{123}, y_{133}, y_{233},\\
  y_{1111}, y_{1112}, y_{1113}, y_{1123}, y_{1133}, y_{1233}, y_{11111}, y_{11112}, y_{11113}, y_{11123}, y_{11133}, y_{11233}\}  $\\
of the $1+3+5+6+6+6=27$ parametric jets up to order $5$.   \\
Using the characters, we get at once $dim(g_0)=1, dim(g_1)= 3, dim(g_2)=5, dim(g_q)=6, \forall q\geq 3$. \\
The localized system $R_4$ with respect to ${\chi}_1$ will have full classes $3$ ($1$ equation) and $2$ ( $4$ equations) of course but the remaining equations will only include the $5$ localized equations of $R_3$. The whole set of ${\alpha}=6$ parametric jets of this new system of finite type will be $\{ y, y_2, y_3, y_{23}, y_{33}, y_{233}\}$. They provide $6$ modular equations generated by the single modular equation $E \equiv a^{233} + {\chi}_1 a^{2223} + ({\chi}_1)^2a^{22222} + ...  =0$ (Exercise).  \\
It is in such a framework that the results of Macaulay look like rather " magic " at first sight, even if they are superseded by the use of formal integrability and involutivity. Indeed, contrary to the first of the examples already presented where new equations of order $2$ could be obtained by differentiating once or even twice the given equations of order $2$ and though the given system is {\it not} homogeneous, the formal integrability of $R_3$ ({\it care}) allows one to use the results of Macaulay on the basis with two generating PD equations of respective orders $3$ and $2$. Accordingly, the formula provided by Macaulay for the Hilbert function in his so-called $H$-basis situation is valid and gives ([9], \S 58, end of p 66) with $n=3$:  \\
\[  \begin{array}{r cl}
 (1-x^3)(1-x^2)(1-x)^{-3} & = & (1+x+x^2)(1+x)(1+x+x^2+x^3+x^4+x^5+...)  \\
   &  =  &  (1+2x+2x^2+x^3)(1+x+x^2+x^3+x^4+x^5+ ...)   \\
   &  =  &  1+3x+5x^2+6x^3+6x^4+6x^5+...  
   \end{array}    \]
in total agreement with the preceding results and the characters obtained.  \\
We finally study the symbols directly. First of all, as $g_4$ is involutive, then it is also $2$-acyclic and we have the long exact sequence: \\
\[  \begin{array}{rcccccccl}
0 \rightarrow & g_6 & \stackrel{\delta}{\longrightarrow} & T^*\otimes g_5 & \stackrel{\delta}{\longrightarrow} & {\wedge}^2T^*\otimes g_4 &  \stackrel{\delta}{\longrightarrow} & {\wedge}^3T^*\otimes g_3 & \rightarrow 0  \\
   0 \rightarrow & 6  &  \rightarrow & 18 & \rightarrow & 18 & \rightarrow & 6 & \rightarrow 0
   \end{array}  \]
However, $g_3$ is not involutive and not even $2$-acyclic because the long exact sequence: \\
\[  \begin{array}{rcccccccl}
0 \rightarrow & g_5 & \stackrel{\delta}{\longrightarrow} & T^*\otimes g_4 & \stackrel{\delta}{\longrightarrow} & {\wedge}^2T^*\otimes g_3 &  \stackrel{\delta}{\longrightarrow} & {\wedge}^3T^*\otimes g_2 & \rightarrow 0  \\
   0 \rightarrow & 6  &  \rightarrow & 18 & \rightarrow & 18 & \rightarrow & 5 & \rightarrow 0
   \end{array}  \]
cannot be exact by counting the dimensions. Indeed, it is known to be exact at $T^* \otimes g_4$ because $g_4$ is the first prolongation of $g_3$ and the image of the central $\delta$ has dimension equal to $18 - 6=12$. Moreover, $g_3$ is strictly contained in the first prolongation of $g_2$ because there is one additional third order equation and the map $\delta$ on the right is well defined. However, as the dimension of the image of the right $\delta$ is smaller or equal to $5$, the dimension of its kernel is thus greater or equal to $18-5=13 > 12$.  \\

\noindent
{\bf EXAMPLE 7}: $n=4 \Rightarrow n+1=5 < 2^4=16=1+4+6+4+1 \Rightarrow dim_k(R_1)=dim_k(J_1(E))=1+4=5,dim_k(R_2)=1+4+6=11,dim_k(R_3)=1+4+6+4=15, dim_k(R_{4+r})=1+4+6+4+1=16, \forall r\geq 0 \Rightarrow g_{5+r}=0, \forall r\geq 0$. In order to fulfill these relations, let us consider the following second order system which is homogeneous and thus formally integrable though not involutive:  \\
\[   {\phi}^1 \equiv y_{44}=0, \hspace{3mm} {\phi}^2 \equiv y_{34}-y_{22}=0, \hspace{3mm} {\phi}^3 \equiv y_{33}=0, \hspace{3mm} {\phi}^4 \equiv y_{24}-y_{11}=0  \]
The $10-4=6$ parametric jets of strict order $2$ are $\{y_{11},y_{12}, y_{13}, y_{14}, y_{22}, y_{23} \}$.\\
The symbol $g_3$ is defined by the $1+3+6+6=16$ equations where we separate the various classes: \\
\[  \begin{array}{l}
y_{444}=0  \\
y_{344}=0, y_{334}=0, y_{333}=0  \\
y_{244}=0, y_{234}-y_{113}=0, y_{233}=0, y_{224}=0, y_{223}=0, y_{222}-y_{113}=0    \\
y_{144}=0, y_{134}-y_{122}=0, y_{133}=0, y_{124}-y_{111}=0, y_{114}=0, y_{112}=0
\end{array}  \]
The $20-16=4$ parametric jets of strict order $3$ are $\{ y_{111}, y_{113}, y_{122}, y_{123}\}$.\\
The symbol $g_4$ is defined by the $1+4+10+19=34=35-1$ equations where we separate again the various classes:  \\
\[ \begin{array}{l}
y_{4444}=0  \\
y_{3444}=0, y_{3344}=0, y_{3334}=0, y_{3333}=0  \\
y_{2444}=0, y_{2344}=0, y_{2334}=0, y_{2333}=0, y_{2244}=0, y_{2234}=0, y_{2233}=0,y_{2224}=0,\\
y_{2223}=0, y_{2222}=0     \\
y_{1444}=0, y_{1344}=0, y_{1334}=0, y_{1333}=0, y_{1244}=0, y_{1234}-y_{1113}=0, y_{1233}=0, \\
y_{1224}=0, y_{1223}=0, y_{1222}-y_{1113}=0, y_{1144}=0, y_{1134}=0, y_{1133}=0, y_{1124}=0,   \\
y_{1123}=0, y_{1122}=0, y_{1114}=0, y_{1112}=0, y_{1111}=0
\end{array}  \]
The only parametric jet of strict order $4$ is $\{y_{1113}\}$.  \\
All the jets of order $\geq 5$ vanish.  \\
We now prove that $g_4$ is $3$-acyclic by showing that the sequence:  \\
\[ \begin{array}{rcccl}
0 \rightarrow & {\wedge}^3T^*\otimes g_4  & \stackrel{\delta}{\longrightarrow} &  {\wedge}^4T^*\otimes g_3   & \rightarrow 0 \\
0 \rightarrow & 4 & \stackrel{\delta}{\longrightarrow} &  4  & \rightarrow 0 
\end{array}  \]
is exact. Computing the central map explicitly, we find:  \\
\[ \begin{array}{rcccl}
 v_{111,1234} & = & v_{1111,234}+v_{1112,341}+v_{1113,412}+v_{1114,123} & = & v_{1113,124}  \\
 v_{122,1234} & = & v_{1122,234}+v_{1222,341}+v_{1223,412}+v_{1224,123} & = & v_{1113,134} \\
 v_{113,1234} & = & v_{1113,234}+v_{1123,341}+v_{1133,412}+v_{1134,123} & = & v_{1113,234} \\
 v_{123,1234} & = & v_{1123,234}+v_{1223,341}+v_{1233,412}+v_{1234,123} & = & v_{1113,123}
 \end{array}  \]
 Accordingly, $g_4$ is $3$-acyclic because $g_{5+r}=0, \forall r\geq 0$ and the central map, which is a monomorphism, is also an isomorphism. However, $g_4$ is not $4$-acyclic, thus involutive, because $g_5=0$ while $dim(g_4)=1$ and the $\delta$-sequence $0 \rightarrow {\wedge}^4T^*\otimes g_4 \rightarrow 0$ cannot therefore be exact. \\
We now prove that $g_4$ is also $2$-acyclic by proving the injectivity of the first $\delta$-map in the short sequence: \\
\[  \begin{array}{rcccccl}
0 \rightarrow & {\wedge}^2T^*\otimes g_4 & \stackrel{\delta}{\longrightarrow} & {\wedge}^3T^* \otimes g_3  & \stackrel{\delta}{\longrightarrow} &{\wedge}^4T^*\otimes g_2 & \rightarrow 0  \\ 
0 \rightarrow & 6 & \stackrel{\delta}{\longrightarrow} & 16 & \stackrel{\delta}{\longrightarrow} & 6 & \rightarrow 0 
\end{array}  \]
Among the $16$ defining equations of the left map $\delta$, we have in particular:  \\
\[  \begin{array}{rcccl}
v_{111,123} & = & v_{1111,23}+v_{1112,31}+v_{1113,12} & = & v_{1113,12}  \\
v_{111,134} & = & v_{1111,34}+v_{1113,41}+v_{1114,13} & = & v_{1113,41} \\
v_{111,234} & = & v_{1112,34}+v_{1113,42}+v_{1114,23} & = & v_{1113,42} \\
v_{113,123} & = & v_{1113,23}+v_{1123,31}+v_{1133,12} & = & v_{1113,23} \\
v_{113,134} & = & v_{1113,34}+v_{1133,41}+v_{1134,13} & = & v_{1113,34} \\
v_{122,123} & = & v_{1122,23}+v_{1222,31}+v_{1223,12} & = & v_{1113,31}
\end{array}   \]
Hence the left $\delta$-map is injective and $g_4$ is $2$-acyclic because $g_5=0$. \\
It follows that $g_3$ is not $3$-acyclic by counting the dimension.\\
The most difficult result will be to prove that $g_3$ is nevertheless also $2$-acyclic by showing the exactness of the long exact sequence:  \\
\[  \begin{array}{rcccccccl}
0 \rightarrow & T^* \otimes g_4 & \stackrel{\delta}{\longrightarrow} & {\wedge}^2T^*\otimes g_3 & \stackrel{\delta}{\longrightarrow} & {\wedge}^3T^*\otimes g_2 & \stackrel{\delta}{\longrightarrow} & {\wedge}^4T^*\otimes T^* & \rightarrow 0  \\
0 \rightarrow & 4 & \stackrel{\delta}{\longrightarrow} & 24 & \stackrel{\delta}{\longrightarrow} & 24 & \stackrel{\delta}{\longrightarrow} & 4 & \rightarrow 0  
\end{array}   \]
The central map $\delta$ is defined by a $24\times 24$ matrix. We obtain in particular:  \\
\[  \begin{array}{rcccl}
v_{11,123} & = & v_{111,23}+v_{112,31}+v_{113,12} & = & v_{111,23}+v_{113,12} \\
v_{11,124} & = & v_{111,24}+v_{112,41}+v_{114,12} & = & v_{111,24}
\end{array}  \]
and so on. Cancelling these equations, we obtain $12$ equations of the second type, namely:  \\
$v_{111,12}=0,v_{111,14}=0, v_{111,24}=0, v_{113,23}=0, v_{113,24}=0, v_{113,34}=0,   \\
v_{122,13}=0, v_{122,14}=0, v_{122,34}=0, v_{123,12}=0, v_{123,13}=0, v_{123,23}=0  $  \\
and $12$ equations of the first type, namely:  \\ 
\[  \begin{array}{c}
v_{122,12}+v_{123,41}=0, v_{122,12}+v_{111,31}=0, v_{123,41}+v_{111,13}=0, \\
 v_{111,23}+v_{113,12}=0, v_{111,23}+v_{123,42}=0, v_{113,12}+v_{123,24}=0, \\
 v_{ 111,34}+v_{113,41}=0, v_{111,34}+v_{122,42}=0, v_{113,41}+v_{122,24}=0, \\
 v_{113,31}+v_{123,34}=0, v_{113,31}+v_{122,23}=0, v_{122,23}+v_{123,34}=0
 \end{array}  \]
 However, we notice that:
 \[   (v_{122,12}+v_{123,41}) - (v_{122,12}+v_{111,31})=v_{123,41}+v_{111,13} \]
 and so on. We obtain therefore $12+(12-4)=20$ linearly independent equations and the kernel of the central $\delta$-map has a dimension equal to $4$ as wished. Accordingly, the last map $\delta$ in the sequence is surjective. Such  result can also be proved by considering the following commutative and exact diagram where 
 $F_0$ is defined by the short exact sequence:\\
  \[  \begin{array}{rccccccl}
0 \rightarrow & g_2 & \rightarrow & S_2T^* & \rightarrow & F_0& \rightarrow 0  \\
0 \rightarrow & 6 & \rightarrow & 10 &\rightarrow & 4& \rightarrow 0  
\end{array}  \]
with $dim(F_0)=4$ while $F_1$ is defined by the long exact sequence: \\
 \[  \begin{array}{rccccccccl}
0 \rightarrow & g_4 & \rightarrow & S_4T^* &\rightarrow & S_2T^*\otimes F_0& \rightarrow &  F_1 & \rightarrow 0  \\
0 \rightarrow & 1 & \rightarrow & 35 &\rightarrow & 40& \rightarrow &  6 & \rightarrow 0  
\end{array}  \]
with $dim(F_1)=40+1-35=6$. Using the snake theorem of homological algebra, we obtain therefore the following commutative and exact diagram:  \\
\[  \begin{array}{rcccccccl} 
  &  &  & 0 &  & 0 & & 0 &  \\
  &  &  & \downarrow \ &  & \downarrow &  & \downarrow &  \\
  & 0  &\rightarrow & S_5 T^*&\rightarrow &   S_3T^*\otimes F_0 & \rightarrow & T^*\otimes F_1 &\rightarrow 0 \\
  & \downarrow & & \hspace{2mm}\downarrow \delta & & \hspace{2mm}\downarrow \delta & & \parallel & \\
  0 \rightarrow & T^* \otimes g_4 & \rightarrow & T^*\otimes S_4T^* & \rightarrow & T^*\otimes S_2T^*\otimes F_0 & \rightarrow & T^*\otimes F_1 & \rightarrow 0 \\
 & \hspace{2mm} \downarrow \delta &  & \hspace{2mm}\downarrow \delta & & \hspace{2mm} \downarrow \delta & & \downarrow &  \\ 
 0 \rightarrow & {\wedge}^2T^*\otimes g_3 & \rightarrow & {\wedge}^2T^*\otimes S_3T^* & \rightarrow & {\wedge}^2T^* \otimes T^* \otimes F_0 & \rightarrow & 0 &  \\
 & \hspace{2mm}\downarrow \delta & &\hspace{2mm} \downarrow \delta & & \hspace{2mm} \downarrow \delta &  &  &  \\
 0 \rightarrow & {\wedge}^3T^* \otimes g_2 & \rightarrow & {\wedge}^3T^*\otimes S_2T^* & \rightarrow & {\wedge}^3T^*\otimes F_0 & \rightarrow & 0 &  \\
  & \hspace{2mm} \downarrow \delta & & \hspace{2mm}\downarrow \delta & & \downarrow & & & \\
  0 \rightarrow & {\wedge}^4T^* \otimes T^* & = & {\wedge}^4T^*\otimes T^* & \rightarrow & 0 & & & \\
   & \downarrow & & \downarrow & & & & &  \\
   &     0  & &  0 & & & & &  
\end{array}     \]
because $g_3$ is $2$-acyclic and we have therefore the upper short exact sequence:  \\
\[ 0 \rightarrow 56 \rightarrow 80 \rightarrow 24 \rightarrow 0 \] 
In actual practice, as the system is homogeneous with constant coefficients, the $4$ second order operators involved commute $2$ by $2$ and the $6$ resulting second order CC are similar to $d_{44}{\phi}^2 - (d_{34}-d_{22}){\phi}^1=0$.  \\
 From a purely algebraic point of view, the corresponding polynomial ideal $\mathfrak{q}=(({\chi}_4)^2,{\chi}_3{\chi}_4-({\chi}_2)^2, ({\chi}_3)^2, {\chi}_2{\chi}_4-({\chi}_1)^2)\subset k[\chi]$ is primary because it is easy to check that $rad(\mathfrak{q})=({\chi}_1,{\chi}_2,{\chi}_3,{\chi}_4)=\mathfrak{m}\in max(k[\chi])$.\\
In order to explain on this example the "translation" of the language of Macaulay into the modern language of differential modules, we may define a differential module $M$ by the exact sequence $  D^6 \rightarrow D^4 \rightarrow D \rightarrow M \rightarrow 0$ over the ring of differential operators $D=k[d_1,d_2,d_3,d_4]=k[d]$ and introduce the corresponding dual system $R=M^*=hom_K(M,k)\simeq hom_D(M,D^*)$. According to the previous results, $dim_k(R)=8< \infty $ and $M$ is $4$-pure. \\ 
  
Finally, having in mind the first of the previous examples, it is not easy to provide a situation where the new system is no longer homogeneous by adding first order terms, but nevertheless formally integrable. For this, let us consider the new second order system $R'_2$:  \\ 
\[    y_{44}=0, \hspace{3mm} y_{34}-y_{22}-y_1=0, \hspace{3mm}  y_{33}=0, \hspace{3mm}  y_{24}-y_{11}-y_3=0  \] 
The $1+3+6+6=16$ third order equations of the first prolongation $R'_3$ are:  \\
\[  \begin{array}{l}
y_{444}=0  \\
y_{344}=0, y_{334}=0, y_{333}=0  \\
y_{244}=0, y_{234} - y_{113}=0, y_{233}=0, y_{224} + y_{14}=0, y_{223} + y_{13}=0, y_{222} - y_{113} + y_{12}=0    \\
y_{144}=0, y_{134} - y_{122} - y_{11}=0, y_{133}=0, y_{124} - y_{111} - y_{13}=0, y_{114} + y_{22} + y_1=0, \\
y_{112} + y_{23} + y_{14}=0
\end{array}  \]
Similarly, the  $1+4+10+19=34$ fourth order equations of $R'_4$ are:  \\
\[ \begin{array}{l}
y_{4444}=0  \\
y_{3444}=0, y_{3344}=0, y_{3334}=0, y_{3333}=0  \\
y_{2444}=0, y_{2344}=0, y_{2334}=0, y_{2333}=0, y_{2244}=0, y_{2234} + y_{122} + y_{11}=0, y_{2233}=0,\\
y_{2224} + y_{111} + y_{13}=0, y_{2223} + y_{123} + y_{11}=0, y_{2222} + 2 y_{122} + y_{11}=0     \\
y_{1444}=0, y_{1344}=0, y_{1334}=0, y_{1333}=0, y_{1244}=0, y_{1234}-y_{1113}=0, y_{1233}=0, \\
y_{1224} - y_{22}-y_1=0, y_{1223} + y_{113}=0, y_{1222}-y_{1113}-y_{23}-y_{14}=0, y_{1144}=0, \\
y_{1134}=0, y_{1133}=0, y_{1124} +y_{113}=0, y_{1123} + y_{122}+ y_{11}=0, y_{1122} + y_{111}=0,\\
 y_{1114} + y_{122} + y_{11}=0, y_{1112} + y_{123} - y_{22} - y_1=0, y_{1111} + 2 y_{113}=0
\end{array}  \] 
The novelty is that, contrary to the previous situation where the second order system $R_2$ considered was {\it automatically} formally integrable because it was homogeneous, {\it now we have to prove the formal integrability of the new system} $R'_2$. However, the symbols $g'_2, g'_3, g'_4$ are easily seen to be identical to the previous corresponding symbols $g_2, g_3, g_4$ while the projections ${\pi}^3_2: R'_3 \rightarrow R'_2$ and ${\pi}^4_3:R'_4 \rightarrow R'_3$ are epimorphisms. Accordingly, as $g'_3=g_3$ is $2$-acyclic, the Janet/Goldschmidt/Spencer criterion of formal integrability proves that $R'_3$ and thus $R'_2$ too are formally integrable, a result leading thus to $g'_5=g_5=0$ as before. \\
It is quite surprising that the type of calculus needed in order to study this purely academic example are {\it exactly} similar to the ones needed in order to study the conformal Killing equations in mathematical physics (See [20],[21],[22]) for more details on this delicate question). \\
 
\noindent
{\bf EXAMPLE 8}: With $n=3$, let us consider the mixed ideal $\mathfrak{a}=({\chi}_1{\chi}_3, {\chi}_2{\chi}_3)=({\chi}_3)\cap ({\chi}_1,{\chi}_2)={\mathfrak{p}}_1 \cap {\mathfrak{p}}_2$. In order to look for the rank of $\mathfrak{a}$ in the sens of Macaulay, it is better to transform this ideal into the system $y_{13}=0, y_{23}=0$ which is {\it not} involutive at first sight and to look for the codimension of the corresponding differential module $M$. Changing $d_1$ to $d_1-d_3$, we get the second order involutive system:  \\
\[ \left\{   \begin{array}{l}
y_{33}-y_{13}=0,  \\
y_{23}=0  
\end{array}
\right. \fbox{ $ \begin{array}{lll}
1 & 2 & 3 \\
1 & 2 & \bullet 
\end{array} $ }      \] 
with full class $3$ only and characters ${\alpha}^3_2=0, {\alpha}^2_2=\alpha=1, {\alpha}^1_2=3$, a result leading to $cd(M)=1$. The localized system $S$ is simply $y_3=0$ with $dim_{k'}(R)= 1$ when $k'=k({\chi}_1,{\chi}_2)$ and $M$ is {\it not} pure because we have the new torsion element ${\bar{y}}_3$ which is killed by $d_2$ and $d_3-d_1$. \\
Setting now $z^1=y, z^2=y_1, z^3=y_2, z^4=y_3$, we get the first order involutive system:  \\
\[ \left\{   \begin{array}{l}
z^1_3-z^4=0, z^2_3-z^3_1=0, z^3_3=0, z^4_3-z^4_1=0,  \\
z^1_2-z^3=0, z^2_2-z^3_1=0, z^4_2=0, \\
z^1_1-z^2=0
\end{array}
\right. \fbox{ $ \begin{array}{lll}
1 & 2 & 3 \\
1 & 2 & \bullet \\
1 & \bullet & \bullet
\end{array} $ }      \] 
with characters ${\alpha}^3_1=0, {\alpha}^2_1=\alpha=1, {\alpha }^1_1=3$. Localizing the equations of class $2$ and $1$ with respect to $k'=k({\chi}_1, {\chi}_2)$, we get:  \\
\[    z^3={\chi}_2 z^1, {\chi}_2z^2-{\chi}_1z^3=0, {\chi}_2 z^4=0 \Rightarrow z^2={\chi}_1z^1 \]
and we notice that {\it the equations of class $1$ are killed by the localization} ([19], Proposition 5.7). Hence, $M$ is not $1$-pure because because we have the new torsion element ${\bar{z}}^4$ which is killed by $d_2$ and $d_3-d_1$. \\

\vspace*{2cm}

\noindent
{\bf BIBLIOGRAPHY}\\
 
\noindent
[1]  ASSEM, I.: Alg\`{e}bres et Modules, Masson, Paris, 1997.\\
\noindent
[2]  BJORK, J.E. : Analytic D-modules and Applications, Kluwer, 1993.\\
\noindent
[3]  BOURBAKI, Alg\`{e}bre Commutative, Chapitre 1 \`{a} 4, Masson, Paris, 1985.\\
\noindent
[4]  BOURBAKI, Alg\`{e}bre, Chapitre 10, Alg\`{e}bre commutative, Masson, Paris, 1980.\\
\noindent
[5]  EISENBUD, D. : Commutative Algebra With a View Towards Algebraic Geometry, Graduate Texts in Math 150, Springer, 1996, (in particular chapter 21).\\
\noindent
[6]  JANET, M.: Sur les Syst\`emes aux d\'eriv\'ees partielles, 
Journal de Math., 8, 3, 1920, 65-151.\\
\noindent
[7]  KASHIWARA, M.: Algebraic Study of Systems of Partial Differential Equations, M\'emoires de la Soci\'et\'e Math\'ematique de France 63, 1995, 
(Transl. from Japanese of his 1970 Master's Thesis).\\
\noindent
[8]  KUNZ, E.: Introduction to Commutative Algebra and Algebraic Geometry, 
Birkh\"{a}user, 1985.\\
\noindent
[9]  MACAULAY, F.S.:  The Algebraic Theory of Modular Systems, Cambridge Tracts 19, Cambridge University Press, London, 1916; Reprinted by Stechert-Hafner Service Agency, New York, 1964.\\
\noindent
[10]  NORTHCOTT, D.G.:  Lessons on Rings, Modules and Multiplicities, Cambridge University Press, 1968.\\
\noindent
[11]  OBERST, U.: Multidimensional Constant Linear Systems, Acta Appl. Math., 20, 1990, 1-175.\\
\noindent
[12]  OBERST, U.:  The Computation of Purity Filtrations over Commutative Noetherian Rings of Operators and their Applications to Behaviours, Multidim. Syst. Sign. Process. (MSSP), Springer, 2013.\\
http://dx.doi.org/10.1007/s11045-013-0253-4  \\
\noindent
[13]  PALAMODOV, V.P.:  Linear Differential Operators with Constant Coefficients,
Grundlehren der Mathematischen Wissenschaften 168, Springer, 1970.\\
\noindent
[14]  POMMARET, J.-F.:  Systems of Partial Differential Equations and Lie Pseudogroups, Gordon and Breach, New York, 1978 
(Russian translation by MIR, Moscow, 1983) \\
\noindent
[15] J.-F. POMMARET, Partial Differential Equations and Group Theory,New Perspectives for Applications, Mathematics and its Applications 293, Kluwer, 1994.\\
http://dx.doi.org/10.1007/978-94-017-2539-2   \\
\noindent
[16]  POMMARET, J.-F.:  Partial Differential Control Theory, Kluwer, 2001, 957 pp.\\
(http://cermics.enpc.fr/$\sim$pommaret/home.html)\\
\noindent
[17]  POMMARET, J.-F.:  Algebraic Analysis of Control Systems Defined by Partial Differential Equations, in Advanced Topics in Control Systems Theory, Lecture Notes in Control and Information Sciences 311, Chapter 5, Springer, 2005, 155-223.\\
\noindent
[18]  POMMARET, J.-F.:  Macaulay Inverse Systems Revisited, Journal of Symbolic Computation, 46, 2011, 1049-1069.  \\
http://dx.doi.org/10.1016/j.jsc.2011.05.007    \\
\noindent
[19]  POMMARET, J.-F.:  Relative Parametrization of Linear Multidimensional Systems, Multidim Syst Sign Process (MSSP), Springer, 2013. \\
http://dx.doi.org/10.1007/s11045-013-0265-0  \\
\noindent
[20]  POMMARET, J.-F.:  The Mathematical Foundations of General Relativity Revisited, Journal of Modern Physics, 2013, 4,223-239.  \\
http://dx.doi.org/10.4236/jmp.2013.48A022  \\
\noindent
[21]  POMMARET, J.-F.:  The Mathematical Foundations of Gauge Theory Revisited, Journal of Modern Physics, 2014, 5, 157-170.  \\
http://dx.doi.org/10.4236/jmp.2014.55026    \\
\noindent
[22]  POMMARET, J.-F.:  Clausius/Cosserat/Maxwell/Weyl Equations: \\
The Virial Theorem Revisited:\\
http://arxiv.org/abs/1504.04118  \\
\noindent
[23]  QUADRAT, A.:  http://wwwb.math.rwth-aachen.de/OreModules  \\
http://www.risc.uni-linz.ac.at/about/conferences/aaca09/   , in particular ... /ModuleTheoryI.pdf   and   ... /ModuleTheoryII.pdf  \\
\noindent
[24]  QUADRAT, A.:  Grade Filtration of Linear Functional Systems, Acta Applicandae Mathematicae, 127, 2013, 27-86. \\
\noindent
[25]  QUADRAT,A.; ROBERTZ, D.:  A Constructive Study of the Module Structure of Rings of Partial Differential Operators, 
Acta Applicandae Mathematicae, 2014, (to appear). \\
http://hal-supelec.archives-ouvertes.fr/hal-00925533  \\
\noindent
[26]  ROTMAN, J.J.: An Introduction to Homological Algebra, Pure and Applied Mathematics, Academic Press, 1979.\\
\noindent
[27]  SCHNEIDERS, J.-P.: An Introduction to D-Modules, Bull. Soc. Roy. Sci. Li\`{e}ge, 63, 1994, 223-295.  \\
\noindent
[28]  SEILER, W.M.: Spencer Cohomology, Differential Equations and Pommaret Bases, Radon Deries Comp. Appl. Math. 1, 2007, 1-50.  \\
\noindent
[29]  SEILER, W.M.: Involution: The Formal Theory of Differential Equations and its Applications to Computer Algebra, Springer, 2009, 660 pp.  \\
 (See also doi:10.3842/SIGMA.2009.092 for a recent presentation of involution, in particular sections 3 (p 6 and reference [11], [22]) and 4). \\
\noindent
[30]  SPENCER, D.C.: Overdetermined Systems of Partial Differential Equations, Bull. Amer. Math. Soc., 75, 1965, 1-114.\\

\end{document}